\documentclass[journal,draft,onecolumn,letterpaper,twoside,11pt]{IEEEtran}

\usepackage{graphicx}
\usepackage{amsfonts,amstext,amsmath,amssymb,amsthm}
\usepackage{bbm}
\usepackage{mathrsfs}

\newcommand{\ind}[1]{\mathbbm{1}_{#1}}   %indicator
\newcommand{\Exp}{{\sf E}}
\newcommand{\Pro}{{\sf P}}
\newcommand{\Real}{\mathbb{R}}

\newcommand{\bH}{{\sf H}}
\newcommand{\cC}{\mathscr{C}}
\newcommand{\ccD}{\mathscr{D}}
\newcommand{\cS}{{\mit \Theta}}
\newcommand{\ccS}{\mathcal{S}}

\newcommand{\cA}{\mathcal{A}}
\newcommand{\cB}{\mathcal{B}}
\newcommand{\cD}{\mathcal{D}}
\newcommand{\cJ}{\mathcal{J}}
\newcommand{\cX}{\mathcal{X}}

\newcommand{\hyptest}{{\renewcommand{\arraystretch}{0.8}\begin{array}{c}\mbox{\raisebox{-0.04em}{\footnotesize $\bH_1$}}\\ \gtreqqless \\ \mbox{\raisebox{0.1em}{\footnotesize $\bH_0$}}
\end{array}}}

\newtheorem{theorem}{Theorem}
\newtheorem{corollary}{Corollary}

\begin{document}

\title{Finite Sample Size Optimality of GLR Tests}

\author{George~V.~Moustakides,~\IEEEmembership{Senior~Member,~IEEE}%
\thanks{G.V. Moustakides is with the Department
of Electrical and Computer Engineering, University of Patras, 26500 Rion,
Greece, e-mail: moustaki@upatras.gr.}%
\thanks{Manuscript received ~~~~, 2009; revised ~~~~, 2009.}}

\markboth{IEEE Transactions on Information Theory,~Vol.~~, No.~~, ~~~2009 (revised)}%
{Moustakides: Finite Sample Size Optimality of GLR Tests}


\maketitle

\begin{abstract}
In several interesting applications one is faced with the problem of simultaneous binary hypothesis testing and parameter estimation.
Although such joint problems are not infrequent, there exist no systematic analysis in the literature that treats them effectively.
Existing approaches consider the detection and the estimation subproblems separately, applying in each case the corresponding optimum
strategy. As it turns out the overall scheme is not necessarily optimum since the criteria used for the two parts are usually incompatible.
In this article we propose a mathematical setup that considers the two problems jointly. Specifically we propose a meaningful combination
of the Neyman-Pearson and the Bayesian criterion and we provide the optimum solution for the joint problem. In the resulting optimum scheme
the two parts interact with each other, producing detection/estimation structures that are completely novel. Notable
side-product of our work is the proof that the well known GLR test is finite-sample-size optimum under this combined sense.
\end{abstract}


\section{Introduction}
\IEEEPARstart{T}{here} exist applications in practice where one must resolve the following problem: decide between two hypotheses $\bH_0$ and $\bH_1$ and then, depending on the decision, estimate a corresponding set of parameters $\theta_0$ or $\theta_1$. Characteristic example of a problem that can be formulated under this combined detection/estimation framework is \textit{target detection and localization} by MIMO radar, where one is not only interested in the classical radar detection problem (presence/absence of a target) but also in estimating its position every time a target is declared present \cite{Bekkerman,Lehmann}. A second example is \textit{retrospective changepoint detection} where we are interested in determining whether there is a point in our samples after which the statistical behavior of the data has changed and, once it is detected then localize this point of interest \cite{Vexler,Gombay}. Clearly segmentation problems can by formulated as retrospective changepoint detection problems.

We would like to emphasize that our goal is not to solve the pure detection problem in the presence of unknown parameters (for this case the parameter estimation subproblem constitutes only an auxiliary step). In our approach the estimation part is a \textit{vital goal} in the whole setup and of the same importance as the detection part. This is clearly the case in the two examples we mentioned before, where the localization of the target in the first and of the changepoint in the second, are of the same importance as the detection part. Current literature does not treat combined problems systematically and the aim of this article is to cover exactly this gap. 

Before introducing in a formal way the combined problem, let us first recall, briefly, the corresponding formulation and the available finite-sample-size optimality results for detection and parameter estimation. For both problems we assume the existence of a random data vector $\cX\in\Real^N$ of length $N$. 

\textit{\underline{Binary hypothesis testing}:} We consider the following two hypotheses $\bH_0,\bH_1$ for $\cX$
\begin{equation}
\bH_i:~~~\cX\sim f_i(X|\theta_i),~i=0,1,
\end{equation}
where ``$\sim$'' means ``distributed according to'' and $f_i(X|\theta_i),~i=0,1,$ are two distinct pdfs with $\theta_i$ denoting a vector of parameters under each hypothesis. Given a realization $X$ of $\cX$, one must decide between the two hypotheses $\bH_0$ and $\bH_1$. If $d\in\{0,1\}$ denotes our decision, then under a Neyman-Pearson formulation we are interested in the following constrained minimization problem
\begin{equation}
\min\Pro(d=0|\bH_1),~\text{subject to}~\Pro(d=1|\bH_0)\le\alpha,
\label{eq:np}
\end{equation}
where $\Pro(\cdot)$ denotes probability and $\alpha\in(0,1)$ the maximal allowable false alarm level.
Optimization is performed over all decision strategies that satisfy the constraint.

Under a {\it finite-sample-size} setting, when the two pdfs are completely known, i.e.~there are no unknown parameters, the optimum test is the celebrated Likelihood Ratio test. If the pdfs have {\it unknown} parameters, except the very rare case where a uniformly most powerful test can be found, the problem in \eqref{eq:np} is not well defined and one needs to resort to min-max formulations for which no systematic solution exists. In this case it is very common to use the Generalized Likelihood Ratio (GLR) test
\begin{equation}
\frac{\sup_{\theta_1\in\cS_1}f_1(X|\theta_1)}{\sup_{\theta_0\in\cS_0}f_0(X|\theta_0)}\hyptest\lambda,
\label{eq:GLRtest}
\end{equation}
where $\mit\Theta_i$ denotes some a-priori known set of values for $\theta_i$.
For the GLR test there is no finite-sample-size optimality result. In fact there are counterexamples against this claim
\cite{Robey,Bose1,Bose2}. Nevertheless the use of the GLR test is widespread in applications and one important sideproduct of our analysis is the demonstration that this popular detection scheme is in fact finite-sample-size optimum {\it under the combined detection and estimation formulation} we are proposing here. We would like to stress that this is no direct contradiction with the counterexamples reported in \cite{Robey,Bose1,Bose2} since in these references the GLR test is evaluated as a pure detector and not in the combined sense we are proposing in this article.

Regarding the problem in \eqref{eq:np}, if we assume that the parameters $\theta_i$ are \textit{random} with known prior pdfs $\pi_i(\theta_i),i=0,1$, then again \eqref{eq:np} has a well defined solution which is the likelihood ratio test between the two marginal pdfs $f_i(X)=\int f_i(X|\theta_i)\pi_i(\theta_i)d\theta_i$.

\textit{\underline{Parameter Estimation}:}
In this problem, we assume that $\cX$ has a pdf $f(X|\theta)$ where $\theta$, as before, denotes a vector of parameters. If $X$ is a realization of $\cX$, the goal is to use the data $X$ in order to provide an estimate $\hat{\theta}$ for $\theta$. Under a finite-sample-size setup, optimum estimation structures are available for the Bayesian formulation and only when $\theta$ is assumed to be random with a known prior pdf $\pi(\theta)$. Specifically, if $C(\hat{\theta},\theta)$ denotes the cost of providing the estimate $\hat{\theta}$ when the true parameter value is $\theta$, then the optimum estimator that \textit{minimizes the average cost} is
\begin{equation}
\hat{\theta}=\text{arg}\inf_U \int C(U,\theta)f(X|\theta)\pi(\theta)d\theta.
\end{equation}
With proper choice of the cost function $C(\hat{\theta},\theta)$, this formula gives rise to a number of well known estimators as the MAP, the conditional mean or the conditional median. 

Next we will combine the two problems and after defining a meaningful performance measure we will develop the optimum detection/estimation structure for the joint problem.

\section{Combined Detection and Estimation}
As we realize from the previous discussion, in both problems, finite-sample-size optimum solutions exist only if we assume that the parameters are random with some known prior. It is therefore natural to expect that the same assumption will be transferred to the more general combined problem. With this observation in mind, let us define the problem of interest.

Consider a random data vector $\cX\in\Real^N$ and the following two hypotheses $\bH_0,\bH_1$:
\begin{equation}
\bH_i:~~\cX\sim f_{i}(X|\theta_i)~\mbox{with
prior pdf}~\pi_{i}(\theta_i),~i=0,1.
\end{equation}
Given any realization $X$ of $\cX$ we would like to decide between the two hypotheses $\bH_0,\bH_1$; and if our decision is in favor of $\bH_i$, then we would like to provide an estimate $\hat{\theta}_i$ for the corresponding parameters $\theta_i$. 

The priors $\pi_i(\theta_i)$ are considered to be generalized functions containing possible point masses. This will allow for the unified analysis of the problem with $\theta_i$ taking a continuum or a discrete set of values. Let us now define what we mean by combined detection/estimation scheme.

\subsection{Combined Detection/Estimation Structure}
We adopt the class of \textit{randomized} detectors and estimators, and we propose the following \textit{two-step scheme}: In the first step with the help of two randomization probabilities $\delta_0(X),\delta_1(X)$ we decide between $\bH_0,\bH_1$. Quantity $\delta_i(X)$ denotes the probability by which we decide $d=i$ using a random game. Clearly $\delta_0(X)+\delta_1(X)=1$.
In the second step we provide parameter estimates that we generate with the help of randomized estimators. Specifically we define two conditional pdfs $q_0(\hat{\theta}_0|X)$ and $q_1(\hat{\theta}_1|X)$, that satisfy
$\int q_0(\hat{\theta}_0|X)\,d\hat{\theta}_0=\int q_1(\hat{\theta}_1|X)\,d\hat{\theta}_1=1$.
These two density functions are applied as follows: if in the first step we decide $d=i$, then in the second step we use the pdf $q_i(\hat{\theta}_i|X)$ to generate a random variable $\hat{\theta}_i$ distributed according to $q_i(\hat{\theta}_i|X)$. This variable constitutes our estimate. Randomized estimators are the direct analog of randomized tests used in hypothesis testing and are not uncommon in Bayesian approaches, as one can verify by consulting \cite[page 65]{Robert}.

We should note that $q_i(\hat{\theta}_i|X)$ must have the same support as the prior $\pi_i(\theta_i)$ since we expect our estimate $\hat{\theta}_i$ to assume the same values as the true parameter $\theta_i$. This is particularly important if $\theta_i$ can take only a finite number of values, in which case $\pi_i(\theta_i)$ and $q_i(\hat{\theta}_i|X)$ will be comprised of point masses. In the latter case, it is easy to see, that we can carry out the analysis using only probabilities instead of pdfs and replace integrals over $\theta_i$ and $\hat{\theta}_i$ with sums.

Summarizing: the combined detection/estimation structure is comprized of the two probabilities $\delta_0(X)$, $\delta_1(X)$ (used in the first step to distinguish between the two hypotheses $\bH_0$, $\bH_1$) and of the two pdfs $q_0(\hat{\theta}_0|X)$, $q_1(\hat{\theta}_1|X)$ (used to provide the necessary parameter estimate in the second step). We denote the complete detection/estimation structure as $\cD=\{\delta_0(X),\delta_1(X),q_0(\hat{\theta}_0|X),q_1(\hat{\theta}_1|X)\}$. 

\textit{Remark 1:} One might wonder if the adoption of a two-step procedure covers all possibilities for a randomized detector/estimator. It turns out that we could also use one-step detectors/estimators that \textit{simultaneously} detect and estimate. However, it is straightforward to show that such schemes can be simulated by properly selected two-step procedures; furthermore, the opposite is also true, that is, any two-step detector/estimator can be simulated by a proper one-step procedure. Consequently the two approaches are fully equivalent and, without loss of generality, we may limit ourselves to the two-step schemes introduced above\footnote{Our claim is particularly easy to prove when the parameters $\theta_i$ take only a finite number of values.}.

In the next subsection our aim is to to define a suitable performance measure for $\cD$ and a corresponding optimization problem that will lead to the identification of the optimum detection/estimation structure.

\subsection{Combined Optimization Problem}
As we mentioned in the Introduction, we are going to combine the Bayesian with the Neyman-Pearson approach.
To this end let $C_{ji}(\hat{\theta}_j,\theta_i)$ denote the cost of deciding in favor of hypothesis $\bH_j$ in the first step and providing the estimate $\hat{\theta}_j$ in the second step, when the true hypothesis is $\bH_i$ and the true parameter is $\theta_i$.

Let us consider the average cost $\cC_i(\cD)$ \textit{given} that the true hypothesis is $\bH_i$. We can express $\cC_i(\cD)$ in terms of the complete detection/estimation structure as follows
\begin{equation}
\cC_i(\cD)=\int\left\{\delta_0(X)\int q_0(\hat{\theta}_0|X)\ccD_{0i}(\hat{\theta}_0,X)d\hat{\theta}_0
+\delta_1(X)\int q_1(\hat{\theta}_1|X)\ccD_{1i}(\hat{\theta}_1,X)d\hat{\theta}_1\right\}dX,
\label{eq:avcost}
\end{equation}
where $\ccD_{ji}(U,X)=\int C_{ji}(U,\theta_i)f_i(X|\theta_i)\pi_i(\theta_i)d\theta_i$. As we can see the four functions $\ccD_{ji}(U,X)$ depend on the known cost functions $C_{ji}(U,\theta_i)$ and on prior information, consequently they are also known and independent from the detection/estimation structure $\cD$. 

We can now define the following optimization problem that we propose as an alternative to the classical problem depicted in \eqref{eq:np}.
\begin{equation}
\inf_{\cD}\cC_1(\cD),~~\text{subject to}~\cC_0(\cD)\le\alpha.
\label{eq:bayesnp}
\end{equation}
Level $\alpha$ constitutes the maximally allowable cost under hypothesis $\bH_0$. As we can see by direct comparison with \eqref{eq:np}, we follow a Neyman-Pearson like approach, having replaced the (conditional) error probabilities of the classical approach with the conditional Bayesian costs.
The problem defined in \eqref{eq:bayesnp} makes a lot of sense. Indeed if one is interested in parameter estimation under each hypothesis then the primal concern is the induced average estimation cost, which quantifies the quality of the corresponding estimate. It is therefore understandable that both, the detection and the estimation subproblems must contribute towards the optimization of the same figure of merit.

Before continuing with the general solution of our problem, we would like to consider a special case which establishes finite-sample-size optimality for the GLR test. The practical significance of this popular test certainly justifies this special analysis. There is however an additional reason that makes this short parenthesis necessary: we plan to use the GLR test as our prototype, therefore we will observe under what conditions we can guarantee its optimality. Then we will apply similar assumptions in the general case, in order to generate GLR-like tests that are compatible with various well known cost functions used in applications. This will produce novel tests that are hopefully more suitable than the classical GLR test, for these problems.

\subsection{Optimality of the GLR Test}
Consider the case where $\theta_i$ takes a finite set of values. Without loss of generality, we will assume that $\theta_i=1,2,\ldots,L_i$ and for simplicity, when $\theta_i=l$, we are going to denote the corresponding pdf as $f_{il}(X)$ instead of $f_i(X|\theta_i=l)$. This immediately suggests that the two prior pdfs $\pi_i(\theta_i)$ will be comprised of an equivalent number of point masses. We denote the corresponding prior probabilities with $\pi_{il}$. In other words under hypothesis $\bH_i$ we have $\cX\sim f_{il}(X)$ with prior probability $\pi_{il}$, where $i=0,1$ and $l=1,\ldots,L_i$. Since $\theta_i$ assumes a finite number of values, the estimators $q_i(\hat{\theta}_i|X)$ will be comprised of point masses as well. Let
$q_{il}(X)$ denote the corresponding probabilities. Our detection/estimation structure can then be identified as the following collection of probabilities
\begin{equation}
\cD=\{\delta_0(X),\delta_1(X),q_{01}(X),\ldots,q_{0L_0}(X),q_{11}(X),\ldots,q_{1L_1}(X)\}
\end{equation}
with the following properties
\begin{equation}
\delta_i(X)\ge0;~q_{il}(X)\ge0;~
\delta_0(X)+\delta_1(X)=\sum_{l=1}^{L_i}q_{il}(X)=1.
\end{equation}
As before the probabilities $\delta_0(X),\delta_1(X)$ are used in the first step to decide between the two main hypotheses. Given that the decision in the first step is in favor of $\bH_i$, we go to the second step and with the help of the probabilities $q_{il}(X),~l=1,\ldots,L_i,$ we decide with the help of a randomized test among the possibilities $f_{i1}(X),\ldots,f_{iL_i}(X)$.

Consider now the following special case of cost functions 
\begin{equation}
C_{10}(\hat{\theta}_1,\theta_0)=C_{01}(\hat{\theta}_0,\theta_1)=1;~
C_{11}(\hat{\theta},\theta)=C_{00}(\hat{\theta},\theta)=\ind{\{\hat{\theta}\ne\theta\}},
\end{equation}
where $\ind{A}$ denotes the indicator of the set $A$. In other words the cost is 0 only when both steps make the correct selection and it is equal to 1 otherwise.
The corresponding average cost $\cC_i(\cD)$ is then equal to the probability of \textit{detection/estimation-error} under hypothesis $\bH_i$. We have the following theorem that solves the problem defined in \eqref{eq:bayesnp}.

\begin{theorem}
Consider the class $\cJ_\alpha$ of all detection/estimation strategies that satisfy the
constraint
\begin{equation}
\Pro(\text{Detection/estimation-error}|\bH_0)\le\alpha,
\label{eq:constr}
\end{equation}
where $\alpha_{\min}<\alpha<1$, with
\begin{equation}
\alpha_{\min}=1-\int\max_{1\le l\le
L_0}\{\pi_{0l}f_{0l}(X)\}dX.
\label{eq:a_min}
\end{equation}
The test, within the class $\cJ_\alpha$, that minimizes the probability
$\Pro(\text{Detection/estimation-error}|\bH_1)$ is given by:

{\it Step\,1}: The optimum strategy for deciding between the two main
hypotheses $\bH_0$ and $\bH_1$ is
\begin{equation}
\frac{\displaystyle\max_{1\le l\le L_1}\{\pi_{1l}f_{1l}(X)\}}{\displaystyle\max_{1\le l\le L_0}\{\pi_{0l}f_{0l}(X)\}}
\hyptest\lambda
\label{eq:glr}
\end{equation}
where, whenever the left hand side coincides with the
threshold we perform a randomization between the two hypotheses and select $\bH_1$ with probability $\gamma$.

{\it Step\,2}: If in Step\,1 we decide in favor of hypothesis $\bH_j$ then
the optimum estimation strategy is
\begin{equation}
\hat{\theta}_j=\text{arg}\max_{1\le l\le L_j}\{\pi_{jl}f_{jl}(X)\}.
\label{eq:opt_iso}
\end{equation}
If more than one indexes attain the same maximum we perform
an arbitrary randomization among them.

The threshold $\lambda$ and the randomization probability $\gamma$
of Step\,1 must be selected so that the constraint in (\ref{eq:constr}) is satisfied
with equality.
\end{theorem}

\begin{IEEEproof} We observe that
$\Pro(\mbox{Detection/estimation-error}|\bH_i)=1-\Pro(\mbox{Correct-detection/estimation}|\bH_i)$,
therefore the constraint is equivalent to
$\Pro(\mbox{Correct-detection/estimation}|\bH_0)\ge1-\alpha$. If we denote the possibility $\{\cX\sim f_{il}(X)\}$ with $\bH_{il}$
then we can write
\begin{equation}
\Pro(\mbox{Correct-detection/estimation}|\bH_i)=\sum_{l=1}^{L_i}\Pro(\mbox{Correct-detection/estimation}|\bH_{il})\pi_{il}
\label{eq:cor_iso_cond}
\end{equation}
with
\begin{equation}
\Pro(\mbox{Correct-detection/estimation}|\bH_{il})=\int\delta_i(X)q_{il}(X)f_{il}(X)dX.
\label{eq:cor_iso}
\end{equation}
Instead of minimizing the probability of detection/estimation-error we can equivalently maximize the probability of correct-detection/estimation. To solve the constrained optimization problem, let
$\lambda>0$ be a Lagrange multiplier and, as in the classical Neyman-Pearson case, with the help of
(\ref{eq:cor_iso_cond}) and (\ref{eq:cor_iso}), define the
corresponding unconstrained version. We then note
\begin{align}
\Pro(\mbox{Correct-detection/estimation}|\bH_1)+\lambda\,\Pro(\mbox{Correct-detection/estimation}|\bH_0)\hskip-10.5cm&\nonumber\\
&=\int\delta_1(X)\left\{
\sum_{l=1}^{L_1}q_{1l}(X)\pi_{1l}f_{1l}(X)
\right\}dX+
\lambda\,\int\delta_0(X)\left\{
\sum_{l=1}^{L_0}q_{0l}(X)\pi_{0l}f_{0l}(X)
\right\}dX\\
&\le\int\delta_1(X)\max_{1\le l\le L_1}\left\{
\pi_{1l}f_{1l}(X)
\right\}dX+
\lambda\,\int\delta_0(X)\max_{1\le l\le L_0}\left\{\pi_{0l}f_{0l}(X)
\right\}dX\label{eq:12}\\
&=\int\left[\delta_1(X)\max_{1\le l\le L_1}\left\{
\pi_{1l}f_{1l}(X)
\right\}+
\delta_0(X)\lambda\max_{1\le l\le L_0}\left\{\pi_{0l}f_{0l}(X)
\right\}\right]dX\\
&\le\int\max\left\{
\max_{1\le l\le L_1}\left\{
\pi_{1l}f_{1l}(X)
\right\},
\lambda\max_{1\le l\le L_0}\left\{\pi_{0l}f_{0l}(X)
\right\}
\right\}dX.\label{eq:14}
\end{align}
Inequality \eqref{eq:12} is valid because the functions
$q_{il}(X),~l=1,\ldots,L_i$ are nonnegative and
complementary (their sum is equal to 1). Inequality \eqref{eq:14} is also true because the same
properties hold for $\delta_i(X),~i=0,1$. Note that the
final expression constitutes an upper bound on the
performance of any detection/estimation rule. Furthermore
this upper bound is attainable by a specific
detection/estimation strategy. Indeed we note that we have
equality in (\ref{eq:12}) when the estimation probabilities
are selected as
\begin{equation}
q_{ik}(X)=\left\{\begin{array}{cl}
1& \mbox{if}~k=\mbox{arg}\min_{1\le l\le
L_i}\{\pi_{il}f_{il}(X)\}\\
0&\mbox{otherwise},
\end{array}
\right.
\end{equation}
and we randomize if there are more than one indexes
attaining the same maximum. This optimum estimation process is the
randomized equivalent of (\ref{eq:opt_iso}). Similarly we have equality in
(\ref{eq:14}) when we select the detection probabilities to
be
\begin{equation}
\delta_1(X)=\left\{\begin{array}{cl}
1& \mbox{if}~\max_{1\le l\le L_1}\left\{
\pi_{1l}f_{1l}(X)
\right\}\ge\lambda\max_{1\le l\le L_0}\left\{\pi_{0l}f_{0l}(X)
\right\}\\
\gamma& \mbox{if}~\max_{1\le l\le L_1}\left\{
\pi_{1l}f_{1l}(X)
\right\}=\lambda\max_{1\le l\le L_0}\left\{\pi_{0l}f_{0l}(X)
\right\}\\
0&\mbox{otherwise},
\end{array}
\right.
\end{equation}
and $\delta_0(X)=1-\delta_1(X)$. Clearly this optimum detection
procedure is the equivalent of (\ref{eq:glr}).

As far as the false alarm constraint is concerned let us
define the following sets
\begin{align}
\begin{split}
\cA(\lambda)&=\left\{X:
\frac{\max_{1\le l\le L_1}\left\{
\pi_{1l}f_{1l}(X)
\right\}}{\max_{1\le l\le L_0}\left\{\pi_{0l}f_{0l}(X)
\right\}}>\lambda
\right\}\\
\cB(\lambda)&=\left\{X:
\frac{\max_{1\le l\le L_1}\left\{
\pi_{1l}f_{1l}(X)
\right\}}{\max_{1\le l\le L_0}\left\{\pi_{0l}f_{0l}(X)
\right\}}=\lambda
\right\}.
\end{split}
\end{align}
For the test introduced above, we can then write that
\begin{align}
\begin{split}
\Pro(\text{Detection/estimation-error}|\bH_0)\hskip-2cm&\\
&=1-\int_{\cA(\lambda)}
\max_{1\le l\le L_0}\left\{\pi_{0l}f_{0l}(X)\right\}dX
-\gamma\int_{\cB(\lambda)}
\max_{1\le l\le L_0}\left\{\pi_{0l}f_{0l}(X)\right\}dX\\
&\ge1-\int_{\cA(\lambda)\cup\cB(\lambda)}
\max_{1\le l\le L_0}\left\{\pi_{0l}f_{0l}(X)\right\}dX\\
&\ge1-\int
\max_{1\le l\le
L_0}\left\{\pi_{0l}f_{0l}(X)\right\}dX=\alpha_{\min}.
\end{split}
\end{align}
The lower bound $\alpha_{\min}$ is clearly attainable in the limit
by selecting $\gamma=1$ and letting $\lambda\to0$.
Also the detection/estimation-error
probability is bounded from above by 1 and we can see that this value
can also be attained in the limit by selecting $\gamma=0$
and letting $\lambda\to\infty$. Existence of a
suitable threshold $\lambda$ and a randomization probability $\gamma$
that assure validity of the false alarm constraint with
equality, as well as, optimality of the resulting test in the
desired sense, can be easily demonstrated following exactly
the same steps as in the classical Neyman-Pearson
case\footnote{In the proof we simply replace the pdfs
$f_i(X)$ with the functions $\max_{1\le l\le
L_i}\{\pi_{il}f_{il}(X)\}$. Even though these functions are not densities, the proof
goes through without change.}. This
concludes the proof.
\end{IEEEproof}

We realize that in order to apply the test in (\ref{eq:glr}) we need knowledge of the
prior probabilities $\pi_{il}$. Whenever this information is
not available we can consider equiprobable subcases
and select $\pi_{il}=1/L_i$.
Under this assumption the optimum test in (\ref{eq:glr}) is
reduced to the familiar form of the GLR test, 
\begin{equation}
\frac{\displaystyle\max_{1\le l\le L_1}f_{1l}(X)}{\displaystyle\max_{1\le l\le L_0}f_{0l}(X)}
\hyptest\lambda,
\label{eq:classicalglr}
\end{equation}
after absorbing the two prior probabilities inside the threshold.

Finally, we should
mention that if hypothesis $\bH_0$ is simple or, if under
hypothesis $\bH_0$ we are not interested in the estimation
problem (therefore we can treat it as simple by forming the
marginal density) then
$\Pro(\mbox{Detection/estimation-error}|\bH_0)$ becomes
the usual false alarm probability with corresponding
$\alpha_{\min}=0$. In other words the false alarm
probability can take any value in the interval $(0,1)$ as in
the classical Neyman-Pearson problem.

\textit{Remark 2:} We observe that the optimum test,
under each main hypothesis, selects the most appropriate subcase with
the help of the MAP selection rule \eqref{eq:opt_iso}. The interesting point is that
this selection is performed independently of the other hypothesis and of the corresponding detection strategy.
This is clearly a very desirable characteristic since it separates
the estimation from the detection problem. In our analysis we are going to provide sufficient conditions that can
guarantee the same property under the general formulation.

\textit{Remark 3:} We obtain the GLR test by assuming that the prior probabilities are \textit{uniform}. We will use the same principle in our general formulation to obtain tests that can be used as alternatives to the classical GLR test.

\section{Optimum Detection/Estimation Scheme}
Let us now continue with the solution of the optimization problem defined in \eqref{eq:bayesnp}. We have the following theorem that provides the desired optimal detection/estimation structure.

\begin{theorem} Consider the class $\cJ_\alpha$ of
detection/estimation structures $\cD$ that satisfy $\cC_0(\cD)\le\alpha$.
The test that minimizes the average cost $\cC_1(\cD)$ within the class $\cJ_\alpha$ is given by
\begin{equation}
\inf_U[\ccD_{01}(U,X)+\lambda \ccD_{00}(U,X)]\hyptest\inf_U[\ccD_{11}(U,X)+\lambda \ccD_{10}(U,X)]
%\ccD_{01}(\hat{\theta}_0,X)-\ccD_{11}(\hat{\theta}_1,X)\hyptest \lambda\left[\ccD_{10}(\hat{\theta}_1,X)-\ccD_{00}(\hat{\theta}_0,X)\right],
\label{eq:th2.1}
\end{equation}
with the optimum estimators defined by
\begin{equation}
\hat{\theta}_j={\rm arg}\inf_U[\ccD_{j1}(U,X)+\lambda\ccD_{j0}(U,X)],~j=0,1,
\label{eq:th2.2}
\end{equation}
and $\lambda>0$ a threshold properly selected to satisfy the corresponding constraint with equality.
\end{theorem}

\begin{IEEEproof} Let $\lambda>0$ be a Lagrange multiplier and consider the unconstraint minimization of the combination $\cC_1(\cD)+\lambda\cC_0(\cD)$. Using \eqref{eq:avcost} we can write
\begin{align}
&\cC_1(\cD)+\lambda\cC_0(\cD)\nonumber\\
&\begin{array}{l}\displaystyle
\hskip1cm=\int\Big\{\delta_0(X)\int q_0(\hat{\theta}_0|X)[\ccD_{01}(\hat{\theta}_0,X)+\lambda \ccD_{00}(\hat{\theta}_0,X)]d\hat{\theta}_0\\
\displaystyle\hskip2.5cm+\delta_1(X)\int q_1(\hat{\theta}_1|X)[\ccD_{11}(\hat{\theta}_1,X)+\lambda \ccD_{10}(\hat{\theta}_1,X)]d\hat{\theta}_1\Big\}dX
\end{array}\\
&\hskip1cm\ge\int\Big\{\delta_0(X)\inf_U[\ccD_{01}(U,X)+\lambda \ccD_{00}(U,X)]+\delta_1(X)\inf_U[\ccD_{11}(U,X)+\lambda \ccD_{10}(U,X)]\Big\}dX\label{eq:th2.pr1}\\
&\hskip1cm\ge\int\min\Big\{\inf_U[\ccD_{01}(U,X)+\lambda \ccD_{00}(U,X)],\inf_U[\ccD_{11}(U,X)+\lambda \ccD_{10}(U,X)]\Big\}dX.\label{eq:th2.pr2}
\end{align}
The inequality in \eqref{eq:th2.pr1} is true because
\begin{align}
\begin{split}
\int q_i(\hat{\theta}_i|X)[\ccD_{i1}(\hat{\theta}_i,X)+\lambda \ccD_{i0}(\hat{\theta}_i,X)]d\hat{\theta}_i&\ge
\inf_U[\ccD_{i1}(U,X)+\lambda \ccD_{i0}(U,X)]\int q_i(\hat{\theta}_i|X)d\hat{\theta}_i\\
&=\inf_U[\ccD_{i1}(U,X)+\lambda \ccD_{i0}(U,X)]
\end{split}
\end{align}
with equality iff $q_i(\hat{\theta}_i|X)$ puts all its probability mass on the choice $\hat{\theta}_i=\text{arg}\inf_U[\ccD_{i1}(U,X)+\lambda \ccD_{i0}(U,X)]$, which is thereby optimum. Similarly we have that \eqref{eq:th2.pr2} is true because $\delta_0(X)+\delta_1(X)=1$, and we have equality iff 
\begin{equation}
\delta_1(X)=\left\{
\begin{array}{cl}
1&\text{if}~\inf_U[\ccD_{01}(U,X)+\lambda \ccD_{00}(U,X)]>\inf_U[\ccD_{11}(U,X)+\lambda \ccD_{10}(U,X)]\\
\gamma&\text{if}~\inf_U[\ccD_{01}(U,X)+\lambda \ccD_{00}(U,X)]=\inf_U[\ccD_{11}(U,X)+\lambda \ccD_{10}(U,X)]\\
0&\text{if}~\inf_U[\ccD_{01}(U,X)+\lambda \ccD_{00}(U,X)]<\inf_U[\ccD_{11}(U,X)+\lambda \ccD_{10}(U,X)],
\end{array}
\right.
\end{equation}
with $0\le\gamma\le1$ and $\delta_0(X)=1-\delta_1(X)$. This is the randomized version of \eqref{eq:th2.1}. This completes the proof.
\end{IEEEproof}

\textit{Remark 4:} For the level $\alpha$ we have $\alpha_{\min}<\alpha<\alpha_{\max}$. It is possible to come up with an expression for $\alpha_{\min}$. Indeed, from \eqref{eq:avcost} it is easy to see that
\begin{align}
\cC_0(\cD)&\ge\int\left\{\delta_0(X)\inf_U\ccD_{00}(U,X)+\delta_1(X)\inf_U\ccD_{10}(U,X)\right\}dX\\
&\ge\int\min\left\{\inf_U\ccD_{00}(U,X),\inf_U\ccD_{10}(U,X)\right\}dX=\alpha_{\min}.
\end{align}
This lower bound is in fact attainable by the \textit{optimum scheme} defined with \eqref{eq:th2.1}, \eqref{eq:th2.2}, if we let $\lambda\to0$. Unfortunately a similar expression for the upper bound $\alpha_{\max}$ was not possible to obtain.

\textit{Remark 5:} As we can see from \eqref{eq:th2.1}, \eqref{eq:th2.2} the optimal solutions for the detection and estimation subproblems are interrelated. If we are interested in the same characteristic we encountered in the GLR test, where the two estimation problems are independent from each other and from the detection part, then the following special form of the cost functions can assure the validity of this property
\begin{equation}
C_{01}(U,\theta_1)=C_{01}(\theta_1)~\text{and}~C_{10}(U,\theta_0)=C_{10}(\theta_0).
\label{eq:c0110}
\end{equation}
Indeed we can see that if \eqref{eq:c0110} is true then $\ccD_{01}(U,X)=\ccD_{01}(X)$ and $\ccD_{10}(U,X)=\ccD_{10}(X)$, which implies that the optimum estimators in \eqref{eq:th2.2} simplify to
\begin{align}
\hat{\theta}_0&={\rm arg}\inf_U[\ccD_{01}(X)+\lambda\ccD_{00}(U,X)]={\rm arg}\inf_U\ccD_{00}(U,X)\\
\hat{\theta}_1&={\rm arg}\inf_U[\ccD_{11}(U,X)+\lambda\ccD_{10}(X)]={\rm arg}\inf_U\ccD_{11}(U,X),
\label{eq:final_estimator}
\end{align}
that is, they coincide with the classical Bayesian estimators which we obtain by treating each estimation problem separately. The optimum detector in
\eqref{eq:th2.1}, under the same assumptions takes the form
\begin{equation}
\ccD_{01}(X)-\inf_U\ccD_{11}(U,X)\hyptest \lambda\left[\ccD_{10}(X)-\inf_U\ccD_{00}(U,X)\right],
\label{eq:final_test}
\end{equation}
which of course relies on the optimum cost values.

\textit{Remark 6:} Observing \eqref{eq:th2.1} and \eqref{eq:th2.2} it seems as if the order of the two steps in our two-step procedure has been reversed. This impression however is not exactly correct. We note that the minimum of a function is unique and it is the two minimal values that are used in \eqref{eq:th2.1}. The actual estimates that realize the two minima, and are depicted in \eqref{eq:th2.2}, are not necessarily unique and therefore we might require randomization which is performed in the second step. But even if the two estimators are deterministic, it is the first step that will dictate which of the two values will be used as our actual parameter estimate. And this selection is performed \textit{after} the detection step. Therefore, strictly speaking, the order is not reversed.

\subsection{Special Case}
We would like now to pay attention to a particular case that is common in applications. Consider under $\bH_1$ that $\cX\sim f_1(X|\theta)$ where $\theta$ a parameter vector with known prior $\pi(\theta)$ and under $\bH_0$ we assume that $\cX\sim f_0(X)$. In other words the pdf under $\bH_0$ is completely known. In fact it is  very common to have $f_0(X)=f_1(X|\theta=0)$. Our goal is to test
$\bH_0$ against $\bH_1$, and whenever we decide in favor of $\bH_1$ to provide an estimate $\hat{\theta}$ for the corresponding parameter vector $\theta$. We should mention that the two application problems discussed in the Introduction, fall under this particular class.

Since parameter estimation is needed only under $\bH_1$,
this suggests that a combined detection/estimation structure will be comprised of the following functions
$\cD=\{\delta_0(X),\delta_1(X),q_1(\hat{\theta}|X)\}$ that satisfy $\delta_j(X)\ge0,~j=0,1$, $q_1(\hat{\theta}|X)\ge0$,
$\delta_0(X)+\delta_1(X)=\int q_1(\hat{\theta}|X)d\hat{\theta}=1$.
The two probabilities $\delta_0(X),\delta_1(X)$ will be used in the first step to decide between the two main hypotheses, while $q_1(\hat{\theta}|X)$ will be employed in the second step to provide the necessary estimate for $\theta$, every time a decision in favor of $\bH_1$ is reached.

Regarding the estimation costs we have the following functions $C_{11}(\hat{\theta},\theta)$, $C_{10}(\hat{\theta})$, $C_{01}(\theta)$ and $C_{00}$. As we can see $C_{00}$ is simply a constant, whereas $C_{10}(\cdot)$ and $C_{01}(\cdot)$ are functions of a single quantity. Consider now the following selection $C_{00}=0$ and $C_{10}(\hat{\theta})=1$, then it is easy to verify that $\cC_0(\cD)=\Pro(d=1|\bH_0)$, i.e.~the probability of false alarm. For this particular selection we have the following interesting corollary of Theorem\,2.

\begin{corollary}
Consider the average cost $\cC_1(\cD)$ under $\bH_1$ defined using the two cost functions $C_{11}(\hat{\theta},\theta)$ and $C_{01}(\theta)$. The optimum detection/estimation structure that minimizes $\cC_1(\cD)$ under the constraint that the false alarm probability $\Pro(d=1|\bH_0)$ is no larger than $\alpha\in(0,1)$, is given by
\begin{equation}
\frac{\ccD_{01}(X)-\inf_U\ccD_{11}(U,X)}{f_0(X)}\hyptest\lambda,
\label{eq:cor2}
\end{equation}
for the optimum detector and
\begin{equation}
\hat{\theta}_1=\text{\rm arg}\inf_U\ccD_{11}(U,X)
\label{eq:cor1}
\end{equation}
for the corresponding optimum estimator. The two functions $\ccD_{11}(U,X)$, $\ccD_{01}(X)$ are defined as follows 
\begin{equation}
\ccD_{11}(U,X)=\int C_{11}(U,\theta)f_1(X|\theta)\pi(\theta)d\theta;~~\ccD_{01}(X)=\int C_{01}(\theta)f_1(X|\theta)\pi(\theta)d\theta.
\end{equation}
\end{corollary}

\subsection{Discussion}
In finite-sample-size optimum detection and estimation the need for the prior pdfs constitutes a very severe weakness. As we mentioned earlier, if this information is not available the corresponding optimization problems must be treated in some min-max context. Unfortunately min-max formulations tend to be very difficult to solve even asymptotically, and no systematic solution exists for the problems of detection and estimation. It is of course clear that the same limitation applies in the case of the more general combined detection/estimation problem. 

A simple (ad-hoc) method to bypass the need for resorting to min-max approaches, is to apply the same idea used to demonstrate optimality for the GLR test, namely assume that the priors are \textit{uniform}. Of course this selection is arbitrary and does not guarantee optimality of the corresponding scheme under any possible min-max sense. On the other hand, it is the only logical choice that reflects our complete lack of knowledge about the priors.
The corresponding tests, examples of which will be seen in the next section, it is expected to have the same weakness as the GLR test, with one major difference: they will be tailored to the specific cost function adopted in the estimation subproblem.

\section{Examples}\label{sec:examples}
In this section we present a number of interesting examples
by selecting various well known forms of cost functions. We
basically concentrate on the popular costs
encountered in the classical Bayesian estimation theory.
We start with the MAP estimate which demonstrates optimality of the GLR test in the continuous case.

\subsection{MAP Detection/Estimation}
Consider the following combination of cost functions
\begin{equation}
C_{01}(U,\theta)=C_{10}(U,\theta)=1;~~
C_{00}(U,\theta)=C_{11}(U,\theta)=\left\{
\begin{array}{cl}
0&\|U-\theta\|\le\Delta\ll 1\\
1&\mbox{otherwise}.
\end{array}
\right.
\end{equation}
We recall from the classical Bayesian estimation theory (see \cite[Page 145]{Poor})
that, as $\Delta\to0$ and assuming sufficient smoothness of
the pdf functions,
the specific selection of costs leads
to the MAP parameter estimation under each main hypothesis. Indeed we observe\footnote{The approximate equality
becomes exact as $\Delta\to0$.}
\begin{equation}
\ccD_{jj}(U,X)\approx\int f_j(X|\theta)\pi_j(\theta)d\theta-f_j(X|U)\pi_j(U)V_j(\Delta)
\end{equation}
where $V_j(\Delta)$ is the volume of a hypersphere of radius $\Delta$ (which can be different for each hypothesis if the two parameter vectors are not of the same length). Substituting in
(\ref{eq:final_test}) yields
\begin{equation}
\frac{\sup_Uf_1(X|U)\pi_1(U)}{\sup_Uf_0(U|X)\pi_0(U)}\hyptest\lambda\frac{V_0(\Delta)}{V_1(\Delta)}=\lambda',
\label{MAP1}
\end{equation}
and the optimum estimator under each hypothesis is the MAP estimator
\begin{equation}
\hat{\theta}_j={\rm arg}\sup_U f_j(X|U)\pi_j(U).
\end{equation}
Similarly for the special case of Corollary\,1 if we define
\begin{equation}
C_{11}(U,\theta)=\left\{
\begin{array}{cl}
0&\|U-\theta\|\le\Delta\ll 1\\
1&\mbox{otherwise},
\end{array}
\right.
\end{equation}
and $C_{01}(\theta)=1$, then $\ccD_{11}(U,X)\approx\int f_1(X|\theta)\pi(\theta)d\theta-f_1(X|U)\pi(U)V_1(\Delta)$ and the optimum test in \eqref{eq:cor2} takes the form
\begin{equation}
\frac{\sup_U f_1(X|U)\pi(U)}{f_0(X)}\hyptest\frac{\lambda}{V_1(\Delta)}=\lambda',
\label{MAP2}
\end{equation}
with the optimum estimator being $\hat{\theta}={\rm arg}\sup_U f(X|U)\pi(U)$.
In both tests (\ref{MAP1}) and (\ref{MAP2}), the threshold $\lambda'$ and the corresponding randomization probability $\gamma$ are selected to satisfy the false alarm constraint with equality. If the prior probabilities $\pi_i(\theta_i),\pi(\theta)$ are unknown and are replaced with the uniform over some prior sets $\mit\Theta_i$ we obtain the classical form of the GLR test depicted in \eqref{eq:GLRtest}.

\subsection{MMSE Detection/Estimation}
Let us now develop the first test that can be used as an
alternative to the GLR test. Consider the following costs
\begin{equation}
C_{01}(U,\theta_1)=C_{01}(\theta_1);~~C_{10}(U,\theta_0)=C_{10}(\theta_0);~~~
C_{00}(U,\theta)=C_{11}(U,\theta)=\|U-\theta\|^2,
\end{equation}
where $C_{01}(\theta_1)$, $C_{10}(\theta_1)$ are functions to be specified in the sequel. Due to the specific form of the costs, the two estimators are independent from each other and also independent from the detection part. Under each main hypothesis the optimum estimator is obtained by minimizing the corresponding mean square error. Consequently the optimum estimator is the conditional mean of the parameter vector given the data vector $X$ (see \cite[Page 143]{Poor}). Specifically we have
\begin{equation}
\hat{\theta}_j=\Exp[\theta_j|X,\bH_j]=\frac{\int\theta_j f_j(X|\theta_j)\pi_j(\theta_j)\,d\theta_j}{\int
f_j(X|\theta_j)\pi_j(\theta_j)\,d\theta_j}.
\end{equation}
The corresponding optimum test after substituting in (\ref{eq:final_test}) takes the form
\begin{equation}
\cA_1(X)\hyptest\lambda\cA_0(X)
\label{MMSE1}
\end{equation}
where
\begin{align}
\begin{split}
\cA_0(X)&=\|\hat{\theta}_0\|^2f_0(X)+\int[C_{10}(\theta_0)-\|\theta_0\|^2]f_0(X|\theta_0)\pi_0(\theta_0)\,d\theta_0\\
\cA_1(X)&=\|\hat{\theta}_1\|^2f_1(X)+\int[C_{01}(\theta_1)-\|\theta_1\|^2]f_1(X|\theta_1)\pi_1(\theta_1)\,d\theta_1\\
f_j(X)&=\int f_j(X|\theta_j)\pi_j(\theta_j)\,d\theta_j.
\end{split}
\end{align}
Selecting $C_{01}(\theta_1)=\|\theta_1\|^2$ and
$C_{10}(\theta_0)=\|\theta_0\|^2$ simplifies the test
considerably yielding
\begin{equation}
\frac{\|\hat{\theta}_1\|^2}{\|\hat{\theta}_0\|^2}\frac{f_1(X)}{f_0(X)}=
\frac{\|\hat{\theta}_1\|^2}{\|\hat{\theta}_0\|^2}
\frac{\int f_1(X|\theta_1)\pi_1(\theta_1)\,d\theta_1}{\int f_0(X|\theta_0)\pi_1(\theta_0)\,d\theta_0}\hyptest\lambda.
\end{equation}
We recognize in the second ratio the likelihood that is used to decide optimally
between the two main hypotheses. By including the first ratio of the two norm square estimates, the test performs simultaneously optimum
detection and estimation.

For the special case of Corollary\,1 it is easy to verify that the corresponding test
takes the form
\begin{equation}
\frac{\|\hat{\theta}_1\|^2f_1(X)+\int[C_{01}(\theta)-\|\theta\|^2]f_1(X|\theta)\pi(\theta)\,d\theta}{f_0(X)}\hyptest\lambda,
\label{MMSE2}
\end{equation}
which, if we select $C_{01}(\theta)=\|\theta\|^2$, simplifies to
\begin{equation}
\|\hat{\theta}\|^2\frac{f_1(X)}{f_0(X)}\hyptest\lambda,
\end{equation}
where $\hat{\theta}=\Exp[\theta|X,\bH_1]=\int\theta f(X|\theta)\pi(\theta)d\theta/\int f(X|\theta)\pi(\theta)d\theta$ and $f_1(X)=\int f(X|\theta)\pi(\theta)d\theta$.

In both tests in (\ref{MMSE1}) and (\ref{MMSE2}), if the priors are not known and are replaced by uniforms, we obtain tests that are
the equivalent of the GLR test for the MMSE criterion.

\subsection{Median Detection/Estimation}
As our final example we present the case of the median estimation where $\theta_i, \hat{\theta}_i,\theta, U$ are scalars and we select the cost functions as follows
\begin{equation}
C_{01}(U,\theta)=C_{01}(\theta);~~C_{10}(U,\theta)=C_{10}(\theta);~~
C_{00}(U,\theta)=C_{11}(U,\theta)=|U-\theta|.
\end{equation}
The estimators are again independent from each other and from detection part. Under each hypothesis we perform optimum Bayes estimation and for this specific cost function
we know that the optimum estimator is the conditional median \cite[Page 143]{Poor}
\begin{equation}
\hat{\theta}_j={\rm \arg}\left\{y:\Pro(\theta_j\le y|X,\bH_j)=\frac{\int_{-\infty}^y f_j(X|\theta_j)\pi_j(\theta_j)\,d\theta_j}{\int
f_j(X|\theta_j)\pi_j(\theta_j)\,d\theta_j}=\frac{1}{2}\right\}.
\end{equation}
The optimum test, as before, becomes
\begin{equation}
\cA_1(X)\hyptest\lambda\cA_0(X)
\end{equation}
where
\begin{align}
\begin{split}
\cA_0(X)&=\int\left[C_{10}(\theta_0)+\theta_0{\rm sgn}(\hat{\theta}_0-\theta_0)\right]f_0(\theta_0|X)\pi_0(\theta_0)d\theta_0\\
\cA_1(X)&=\int\left[C_{01}(\theta_1)+\theta_1{\rm sgn}(\hat{\theta}_1-\theta_1)\right]f_1(\theta_1|X)\pi_1(\theta_1)d\theta_1.
\end{split}
\end{align}
If additionally we select $C_{01}(\theta_1)=|\theta_1|$ and $C_{10}(\theta_0)=|\theta_0|$ then the optimum test takes the more convenient form
\begin{equation}
\frac{\int_0^{\hat{\theta}_1}\theta_1f_1(X|\theta_1)\pi_1(\theta_1)d\theta_1}
{\int_0^{\hat{\theta}_0}\theta_0f_0(X|\theta_0)\pi_0(\theta_0)d\theta_0}\hyptest\lambda.
\end{equation}
For the special case of Corollary\,1 and for $C_{01}(\theta)=|\theta|$, the corresponding optimum test reduces to
\begin{equation}
\frac{\int_0^{\hat{\theta}}\theta f(X|\theta)\pi(\theta)\,d\theta}{f_0(X)}
\hyptest\lambda,
\end{equation}
while the optimum estimator is $\hat{\theta}={\rm arg}\{y:\Pro(\theta\le y|X,\bH_1)=0.5\}$.
Finally when the priors are selected to be uniform, we obtain a test that is the alternative to the GLR test but tuned
to the specific Bayesian criterion we employ in the estimation part.

\section{Application to Retrospective Changepoint Detection}
Perhaps the most appropriate application where one would readily need to replace the GLR test with an alternative scheme, is the problem of target detection and localization. Clearly for this problem the most suitable cost function is the mean square error between the location estimate and the true position. This choice will inevitably lead to the use of tests that are similar to \eqref{MMSE2}, proposing a completely novel approach for this intriguing problem. Unfortunately the corresponding derivations are lengthy and thus impossible to detail here. In the limiting space we have to our disposal it is feasible to treat, with our preceding methodology, the second application we mentioned in the Introduction, namely the retrospective changepoint detection problem. We would like to mention that even though in this problem the estimation costs are MAP-like, suggesting use of the GLR test, as we will see, there is sufficient simplicity and originality in our results that make our analysis interesting and worth including in this article.

In its \textit{simplest} form, retrospective changepoint detection is about an observation vector $\cX\in\Real^N$ and two pdfs $f_\infty(X)$ and $f_0(X)$ which are completely known. If $\cX=[\chi_1,\ldots,\chi_N]^t$ then we assume that there is an unknown point $\tau$ such that the samples $\{\chi_1,\ldots,\chi_\tau\}$ follow the nominal measure $f_\infty(X)$ while the $\{\chi_{\tau+1},\ldots,\chi_N\}$ switch to the alternative $f_0(X)$. Consequently, the changepoint $\tau$ is the \textit{last} point where the samples follow the nominal regime\footnote{Notation seems to be somewhat awkward compared to the usual one used in hypothesis testing. We simply follow the standard practice of sequential changepoint detection theory.}. 

We are interested in deciding whether the change took place within or before the given collection of samples, that is $\tau<N$, or the change will take place at some future point (possibly at infinity), that is $\tau\ge N$. In the former case we would also like to obtain an estimate $\hat{\tau}$ of the changepoint $\tau$. The combined detection/estimation version of the retrospective changepoint detection problem, as we mentioned in the Introduction, is suitable for formulating segmentation problems. 

Let us first define the joint pdf $f_\tau(X)$ of the samples $\cX$ given $\tau$. We distinguish three sets of values for $\tau$, namely $\tau\le0$, $\tau\in\{1,\ldots,N-1\}$ and $\tau\ge N$. The first corresponds to a change occurring before taking any samples, the second to a change within the available sample set and the third to the change occurring after we acquired the samples.
The most common model for the induced joint pdf is \cite{Moustakides1,Moustakides2}
\begin{equation}
f_\tau(X)=\left\{
\begin{array}{cl}
f_0(X)&\text{for}~\tau\le0\\
f_\infty(X_1^\tau)f_0(X_{\tau+1}^N|X_1^n)&\text{for}~0<\tau\le N-1\\
f_\infty(X)&\text{for}~N\le \tau,
\end{array}
\right.
\label{eq:fn}
\end{equation}
where $X=[x_1,\ldots,x_N]^t$ and for $a\le b$ we define $X_a^b=[x_a,\ldots,x_b]^t$.
We can see that if the change takes place before the samples are acquired, all samples are under the alternative regime. If the change takes place within the available set, then the initial portion of the samples follows the pdf of the nominal regime while the final portion the \textit{conditional} pdf of the alternative regime. Finally if the change does not occur before or inside the available data set, all samples are under the nominal regime.

Regarding the changepoint $\tau$ there are different models. Detailed discussion of the various possibilities can be found in \cite{Moustakides1,Moustakides2}. Here we limit ourselves to Shiryaev's popular Bayesian model. Specifically we assume that $\tau$ is a random variable with a prior $\{\varpi_n\}$ defined as $\varpi_0=\Pro(\tau\le0)$, $\varpi_n=\Pro(\tau=n)$ for $0<n\le N-1$, $\varpi_N=\Pro(\tau\ge N)$ and such that $\sum_{n=0}^N\varpi_n=1$. 

As we mentioned, the goal is to test $\{\tau\le N-1\}$ against $\{\tau\ge N\}$, and in the former case provide and estimate $\hat{\tau}$ for $\tau$.
Formulating the problem according to our previous theory, we have that under $\bH_0$ the samples follow the nominal pdf $f_\infty(X)$ while under $\bH_1$ we have $N$ different possibilities with corresponding pdf $f_\tau(X)$ and prior $\pi_\tau=\varpi_\tau/(\sum_{k=0}^{N-1}\varpi_k)=\varpi_\tau/(1-\varpi_N)$, where $0\le \tau< N$.

Let us now consider the combined detection/estimation problem in the sense of Corollary\,1, namely minimize the average cost under $\bH_1$ subject to a false alarm probability constraint under $\bH_0$. We propose the following cost functions $C_{11}(\hat{\tau},\tau)=\ind{\{\hat{\tau}\neq\tau\}}$, where $\ind{A}$ denotes the indicator function of the set $A$, and $C_{01}(\tau)=1$. In other words we penalize with 1 the incorrect detection of $\bH_1$ but also the correct detection of $\bH_1$ followed by an incorrect estimation of $\tau$. The average cost is simply the probability of detection/estimation-error introduced in Subsection\,II.C. Applying the results of Corollary\,1 and using the Bayes rule and \eqref{eq:fn}, the optimum detection/estimation structure is given by
\begin{equation}
\max_{0\le n< N}\left\{\pi_n\frac{f_n(X)}{f_\infty(X)}\right\}=\max_{0\le n< N}\left\{\pi_n\frac{f_0(X_{n+1}^N|X_1^n)}{f_\infty(X_{n+1}^N|X_1^n)}\right\}\hyptest\lambda
\end{equation}
for the optimum detector and
\begin{equation}
\hat{\tau}=\text{arg}\max_{0\le n< N}\left\{\pi_n\frac{f_0(X_{n+1}^N|X_1^n)}{f_\infty(X_{n+1}^N|X_1^n)}\right\}
\end{equation}
for the corresponding optimum estimator.
If the priors $\pi_n,~n=0,\ldots,N-1,$ are unknown and we select them to be equal, then we obtain the GLR test version of the problem
\begin{equation}
\ccS_N=\max_{0\le n< N}\left\{\frac{f_0(X_{n+1}^N|X_1^n)}{f_\infty(X_{n+1}^N|X_1^n)}\right\}\hyptest\lambda
\label{eq:cusum}
\end{equation}
where $\ccS_N$ is known as the \textit{CUSUM statistic} for point $N$. 

The previous result was of course expected since we followed the same formulation as the one used in Subsection\,II.C to prove optimality of the GLR test. Interestingly, our theory allows for the development of alternative detection/estimation structures in a simple and straightforward manner. For example one might argue that the cost $C_{11}(\hat{\tau},\tau)=\ind{\{\hat{\tau}\ne\tau\}}$ is overly stringent and propose as alternative the function
$C_{11}(\hat{\tau},\tau)=\ind{\{|\hat{\tau}-\tau|>m\}}$ where $0\le m\ll N$ is a nonnegative integer. In other words we tolerate errors in the estimate of $\tau$ that do not exceed $m$ points. If $m=0$ the problem is reduced to the case already discussed. Clearly most practical segmentation problems would allow $m>0$.

Again we adopt the setup proposed in Corollary\,1. It is then easy to verify that we obtain the following optimum structure
\begin{equation}
\max_{m\le n< N-m}\left\{\sum_{k=-m}^m\pi_{n+k}\frac{f_{n+k}(X)}{f_\infty(X)}\right\}=\max_{m\le n< N-m}\left\{\sum_{k=-m}^m\pi_{n+k}\frac{f_0(X_{n+k+1}^N|X_1^{n+k})}{f_\infty(X_{n+k+1}^N|X_1^{n+k})}\right\}\hyptest\lambda
\end{equation}
for the detector and
\begin{equation}
\hat{\tau}=\text{arg}\max_{m\le n< N-m}\left\{\sum_{k=-m}^m\pi_{n+k}\frac{f_0(X_{n+k+1}^N|X_1^{n+k})}{f_\infty(X_{n+k+1}^N|X_1^{n+k})}\right\}
\end{equation}
for the estimator. Finally assuming uniform priors for the case where the probabilities $\{\pi_0,\ldots,\pi_{N-1}\}$ are unknown, leads to the test
\begin{equation}
\bar{\ccS}_N=\max_{m\le n< N-m}\left\{\sum_{k=-m}^m\frac{f_0(X_{n+k+1}^N|X_1^{n+k})}{f_\infty(X_{n+k+1}^N|X_1^{n+k})}\right\}\hyptest\lambda,
\end{equation}
which is completely novel and replaces the GLR test in \eqref{eq:cusum}, with $\bar{\ccS}_N$ being clearly different than the CUSUM statistic.

\section{Conclusion}
By introducing a joint detection/estimation formulation
that properly combines the Neyman-Pearson methodology (for detection) and the Bayesian methodology (for estimation), we derived optimum schemes for problems that require simultaneous detection and estimation. 
Important side-product of our analysis is the demonstration that the well known GLR test is finite-sample-size optimum under this joint-problem sense. Furthermore we were able to provide completely novel GLR-type tests, that were derived by replacing the MAP estimation cost function with other well known choices as the mean square or mean absolute estimation error.
Finally, we used our proposed methodology to analyze the problem of retrospective changepoint detection. This led to the development of a novel detection/estimation structure that can replace the CUSUM approach which is obtained when we apply the GLR test.

\section*{Acknowledgment}
The author would like to thank his good friend, Prof. Igor Nikiforov from the
Universit\'e de Technologie de Troyes (UTT), France, for
enlightening discussions.

\end{document}